\documentclass[12pt,twoside]{article}
\begin{document}

\date{}

%\emph{Int. J. Open Problems Compt. Math., Vol. 3, No. 1, March 2010}\

%\emph{ISSN 1998-6262; Copyright \copyright ICSRS Publication, 2010}\

%\emph{www.i-csrs.org}

\centerline{}

\centerline{}

\centerline {\Large{\bf On Fixed-point theorems in Intuitionistic  }}

\centerline{}

\centerline {\Large{\bf Fuzzy metric Space}}

\centerline{}

%% My definition

\centerline{\bf {T.K. Samanta , Sumit Mohinta  and Iqbal H. Jebril}}

\centerline{}

\centerline{Department of Mathematics, Uluberia College, India-711315.} %
\centerline{e-mail: mumpu$_{-}$tapas5@yahoo.co.in}

\centerline{e-mail: sumit.mohinta@yahoo.com}

\centerline{Department of Mathematics, King Faisal University,Saudi Arabia.} %
\centerline{e-mail: iqbal501@hotmail.com}

\begin{abstract}
\textbf{\emph{In this paper, first we have established two sets of sufficient conditions for a mapping to have unique fixed point in a intuitionistic fuzzy metric space and then we have redefined the contraction mapping in a intuitionistic fuzzy metric space and thereafter we proved the Banach Fixed Point theorem.}}
\end{abstract}

\centerline{}

\newtheorem{Theorem}{\quad Theorem}[section]

\newtheorem{Definition}[Theorem]{\quad Definition}

\newtheorem{Corollary}[Theorem]{\quad Corollary}

\newtheorem{Lemma}[Theorem]{\quad Lemma}

\newtheorem{Note}[Theorem]{\quad Note}

\newtheorem{Remark}[Theorem]{\quad Remark}

\newtheorem{Result}[Theorem]{\quad Result}

\newtheorem{Proposition}[Theorem]{\quad Proposition}

\newtheorem{Example}[Theorem]{\quad Example}

\textbf{Keywords:} \emph{Fuzzy Sets, Intuitionistic Fuzzy Set, Intuitionistic Fuzzy Metric, Contraction Mapping, Contractive Sequence, Cauchy sequence. }\newline
\textbf{2010 Mathematics Subject Classification:} 03F55, 46S40.

%===============================================

\section{Introduction}

%===============================================
Fuzzy set theory was first introduce by Zadeh\cite{zadeh} in 1965 to describe the situation in which data are imprecise or vague or uncertain. Thereafter the concept of fuzzy set was generalized as intuitionistic fuzzy set by K. Atanassov\cite{Atanasov} in 1984. It has a wide range of application in the field of population dynamics , chaos control , computer programming , medicine , etc.\\ Using the idea of intuitionistic fuzzy set, Park\cite{Park} introduced the notion of intuitionistic fuzzy metric spaces with the help of continuous t-norms and continuous t-conorms, which is a generalization of fuzzy metric space due to George and Veeramani\cite{Veeramani}. \\ In intuitionistic fuzzy metric space, Mohamad\cite{Abdul} proved Banach Fixed Point theorem. But in his paper\cite{Abdul}, the definition of contractive mapping is not natural and to prove the iterative sequence is a Cauchy sequence, he first proved that every iterative sequence is a contractive sequence and then assumed that every contractive sequences are Cauchy.
\\ In our paper, we have redefined the notion of contraction mapping in a intuitionistic fuzzy metric space and then directly, it has been proved that the every iterative sequence is a Cauchy sequence, that is, we don't need to assume that every contractive sequences are Cauchy sequences. Thereafter we have established the Banach Fixed Point theorem. In this paper, also, we have established another two sets of sufficient conditions for a mapping to have unique fixed point.

%================================================

\section{Preliminaries}

%================================================
We quote some definitions and statements of a few theorems which will be
needed in the sequel.

\begin{Definition} \cite{Schweizer}.
A binary operation \, $\ast \; : \; [\,0 \; , \; 1\,] \; \times \;
[\,0 \; , \; 1\,] \;\, \longrightarrow \;\, [\,0 \; , \; 1\,]$ \, is
continuous \, $t$ - norm if \,$\ast$\, satisfies the
following conditions \, $:$ \\
$(\,i\,)$ \hspace{0.3cm} $\ast$ \, is commutative and associative ,
\\ $(\,ii\,)$ \hspace{0.2cm} $\ast$ \, is continuous , \\
$(\,iii\,)$ \hspace{0.1cm} $a \;\ast\;1 \;\,=\;\, a \hspace{1.2cm}
\forall \;\; a \;\; \varepsilon \;\; [\,0 \;,\; 1\,]$ , \\
$(\,iv\,)$ \hspace{0.1cm} $a \;\ast\; b \;\, \leq \;\, c \;\ast\; d$
\, whenever \, $a \;\leq\; c$  ,  $b \;\leq\; d$  and  $a \, , \, b
\, , \, c \, , \, d \;\, \varepsilon \;\;[\,0 \;,\; 1\,]$.
\end{Definition}

\begin{Definition}
\cite{Schweizer}. A binary operation \, $\diamond \; : \; [\,0 \; ,
\; 1\,] \; \times \; [\,0 \; , \; 1\,] \;\, \longrightarrow \;\,
[\,0 \; , \; 1\,]$ \, is continuous \, $t$-conorm if \,$\diamond$\,
satisfies the
following conditions \, $:$ \\
$(\,i\,)\;\;$ \hspace{0.1cm} $\diamond$ \, is commutative and
associative ,
\\ $(\,ii\,)\;$ \hspace{0.1cm} $\diamond$ \, is continuous , \\
$(\,iii\,)$ \hspace{0.1cm} $a \;\diamond\;0 \;\,=\;\, a
\hspace{1.2cm}
\forall \;\; a \;\; \in\;\; [\,0 \;,\; 1\,]$ , \\
$(\,iv\,)$ \hspace{0.1cm} $a \;\diamond\; b \;\, \leq \;\, c
\;\diamond\; d$ \, whenever \, $a \;\leq\; c$  ,  $b \;\leq\; d$
 and  $a \, , \, b \, , \, c \, , \, d \;\; \in\;\;[\,0
\;,\; 1\,]$.
\end{Definition}

\begin{Result}
\cite{klement}. $(\,a\,)\;$  For any \, $r_{\,1} \; , \;
r_{\,2} \;\; \in\;\; (\,0 \;,\; 1\,)$ \, with \, $r_{\,1} \;>\;
r_{\,2}$, there exist $\;r_{\,3} \; , \; r_{\,4} \;\; \in \;\; (\,0
\;,\; 1\,)$ \, such that \, $r_{\,1} \;\ast\; r_{\;3} \;>\; r_{\,2}$
\, and \, $r_{\,1} \;>\; r_{\,4} \;\diamond\; r_{\,2}.$
\\ $(\,b\,)$ \, For any \, $r_{\,5} \;\,
\in\;\, (\,0 \;,\; 1\,)$ , there exist \, $r_{\,6} \; , \; r_{\,7}
\;\, \in\;\, (\,0 \;,\; 1\,)$ \, such that \, $r_{\,6} \;\ast\;
r_{\,6} \;\geq\; r_{\,5}$ \,and\, $r_{\,7} \;\diamond\; r_{\,7}
\;\leq\; r_{\,5}.$
\end{Result}

\begin{Definition}
\cite{Park}
 Let \,$\ast$\, be a continuous \,$t$-norm ,
\,$\diamond$\, be a continuous \,$t$-conorm and $X$ be any non-empty set.
An \textbf{intuitionistic fuzzy metric} or in short $\textit{\textbf{IFM}}$ on
\,$X$\, is an object of the form \\ $A \;\,=\;\, \{\; (\,(\,x \;,\,y \;,\;
t\,) \;,\; \mu\,(\,x \;,\,y\;,\; t\,) \;,\; \nu\,(\,x \;,\,y \;,\; t\,) \;) \;\,
: \;\, (\,x \;,\,y \;,\; t\,) \;\,\in\;\, X^{2}\times(0\,,\,\infty)\}$\, where $\mu \,,\, \nu\;$ are fuzzy sets on \, $X^{2}\times(0\,,\,\infty)$ , \,$\mu$\, denotes the degree of nearness and
\,$\nu$\, denotes the degree of non$-$nearness of $x$ and $y$ relative
to $t$ satisfying the following conditions $:$ \, for all $x , y , z
\, \in \,X , \, s , t \, > \, 0$ \\
$(\,i\,)$ \hspace{0.10cm}  $\mu\,(\,x \;,y\;, t\,) \;+\; \nu\,(\,x
\;,y\;, t\,) \;\,\leq\;\, 1 \hspace{1.2cm} \forall \;\; (\,x \;,\,y \;,\;
t\,)
\;\,\in\;\,X^{2}\times(0,\infty) ;$
\\$(\,ii\,)$ \hspace{0.10cm}$\mu\,(\,x \;,\,y \;,\; t\,) \;\,>\;\, 0 \, ;$ \\
$(\,iii\,)$ $\mu\,(\,x \;,\,y\;,\; t\,) \;\,=\;\, 1$ \, if
and only if \, $x \;=\;y \,$\,
\\$(\,iv\,)$ $\mu\,(\,x \;,\,y\;,\; t\,) \;\,=\;\,\mu\,(\,y \;,\,x\;,\; t\,);$
\\ $(\,v\,)$
\hspace{0.10cm} $\mu\,(\,x \;,\,y \;,\; s\,) \;\ast\; \mu\,(\,y \;,\,z \;,\; t\,)
\;\,\leq\;\, \mu\,(\,x \;,\,z \;,\; s\;+\;t\,
\,) \, ;$
\\$(\,vi\,)$ $\mu\,(\,x \,,\,y\,,\,\cdot\,) :(0 \,,\;\infty\,)\,\rightarrow
\,(0 \,,\;1]$ \, is continuous;
\\$(\,vii\,)$ \hspace{0.10cm}$\nu\,(\,x \;,\,y \;,\; t\,) \;\,>\;\, 0 \, ;$
\\$(\,viii\,)$ $\nu\,(\,x \;,\,y\;,\; t\,) \;\,=\;\, 0$ \, if
and only if \, $x \;=\;y \,;$\,
\\$(\,ix\,)$ $\nu\,(\,x \;,\,y\;,\; t\,) \;\,=\;\,\nu\,(\,y \;,\,x\;,\; t\,);$
\\ $(\,x\,)$\hspace{0.10cm} $\nu\,(\,x \;,\,y \;,\; s\,) \;\diamond\; \mu\,(\,y \;,\,z \;,\; t\,)
\;\,\geq\;\, \nu\,(\,x \;,\,z \;,\; s\;+\;t\,
\,) \, ;$
\\$(\,xi\,)$ $\nu\,(\,x \,,\,y\,,\,\cdot\,) :(0\,,\;\infty\,)\,\rightarrow
\,(0\,,\;1]$ \, is continuous.
\\\\ If \,$A$\, is a $\textit{\textbf{IFM}}$ on
\,$X$\, , the pair \,$(\,X \,,\, A\,)$\, will be called a \textbf{intuitionistic fuzzy metric space} or in short $\textit{\textbf{IFMS}}$.\\\\ We further assume that \,$(\,X \,,\, A\,)$\, is a $\textit{\textbf{IFMS}}$ with the property \\
$(\,xii\,)\;$ For all \,$a\,\in\,(0 \,,\, 1)$,\, $a \,\ast\, a \,=\, a$\, and \,$a \,\diamond\, a \,=\, a$
\end{Definition}

\begin{Remark}
\cite{Park}
In intuitionistic fuzzy metric space $X$, \, $\mu\,(\,x \,,\,y \,,\,\cdot\,)$ \,
is non-decreasing and \,$\nu\,(\,x \,,\,y \,,\,\cdot\,)$\, is non-increasing for all \,$x\,,\,y\,\in\,X$.
\end{Remark}

\begin{Definition}
\cite{Veeramani}
 A sequence $\{\,x_{\,n}\,\}_{\,n}$ in a intuitionistic fuzzy metric space is said to be a \textit{\textbf{Cauchy sequence}} if and only if for each \,$r\,\in\,(\,0\,,\,1\,)$\, and \,$t>0$\, there exists \,$n _{\,0}\,\in\,N$\, such that \,$\mu\,(\,x_{n} \,,\,x_{m} \,,\, t\,) \;>\;1\,-\,r $\, and \,$\nu\,(\,x_{n} \,,\,x_{m} \,,\, t\,) \;<\;r$\, for all $\;n\,,\,m\;\geq\;n_{\,0}.$\, \\A sequence $\{\,x_{\,n}\,\}$ in a intuitionistic fuzzy metric space is said to converge to \,$x\,\in\,X$\, if and only if for each \,$r\,\in\,(\,0\,,\,1\,)$\, and \,$t>0$\, there exists \,$n _{\,0}\,\in\,N$\, such that \,$\mu\,(\,x_{n} \,,\,x \,,\, t\,) \;>\;1\,-\,r $\, and \,$\nu\,(\,x_{n} \,,\,x \,,\, t\,) \;<\;r$\, for all $\;n\,,\,m\;\geq\;n_{\,0}.$\,
\end{Definition}

\begin{Note}
\cite{Samanta}
A sequence $\{\,x_{\,n}\,\}_{\,n}$ in an intuitionistic fuzzy metric space is a Cauchy sequence if and only if \[\mathop {\lim }\limits_{n\,\, \to \,\,\infty } \,\mu\,(\,x_{n} \,,\,x_{n\,+\,p} \,,\, t\,)\;=\;1\;\; and \;\; \mathop {\lim }\limits_{n\,\, \to \,\,\infty } \,\nu\,(\,x_{n} \,,\,x_{n\,+\,p} \,,\, t\,)\;=\;0\;\]
A sequence $\{\,x_{\,n}\,\}_{\,n}$ in an intuitionistic fuzzy metric space converges
\,$x\,\in\,X$\, if and only if \[\mathop {\lim }\limits_{n\,\, \to \,\,\infty } \,\mu\,(\,x_{n} \,,\,x \,,\, t\,)\;=\;1\;\; and \;\; \mathop {\lim }\limits_{n\,\, \to \,\,\infty } \,\nu\,(\,x_{n} \,,\,x \,,\, t\,)\;=\;0\;\]
\end{Note}

\begin{Definition}
\cite{Abdul}
An intuitionistic fuzzy metric space $(X,A)$ has the property $\spadesuit$ if \\$\mathop {\lim }\limits_{t\;\, \to \;\,\infty } \;\mu\,(\,x \;,\,y \;,\; t\,)\;=\;1\;$ and $\mathop {\lim }\limits_{t\;\, \to \;\,\infty } \;\nu\,(\,x \;,\,y \;,\; t\,)\;=\;0\;$ for all $\;x,y\;\in\;X$
\end{Definition}
%\smallskip

\begin{Definition}
\cite{Abdul}
 Let $(\,X\,,\,A\,)$ be a intuitionistic fuzzy metric space. We will say the mapping \,$f : X \rightarrow X$\, is \textbf{t-uniformly continuous} if for each \,$\varepsilon,$\, with \,$0\, < \,\varepsilon \,< \,1,$\, there exists  \,$0 \,<\, r \,<\, 1,$\, such that \,$\mu\,(\,x \,,\, y \,,\, t\,) \,\geq\,1\,-\,r\, $\, and \,$\nu\,(\,x \,,\, y \,,\, t\,) \,\leq\,r\, $ implies $\mu\,\left(\,f(x) \,,\, f(y) \,,\, t\,\right) \,\geq\,1\,-\,\varepsilon \, $ and \,$\nu\,\left(\,f(x) \,,\, f(y) \,,\, t\,\right) \,\leq\,\varepsilon \, $ for each \,$x , y\,\in\,X$\, and \,$t \,> \,0.$
\end{Definition}

%\smallskip

\begin{Definition}
\cite{Abdul}
Let \,$(\,X \,,\, A\,)$\, be a intuitionistic fuzzy metric space. A mapping \,$f : X \rightarrow X$\, is \textbf{intuitionistic fuzzy contractive} if there exists \,$k\in\,(0 \,,\, 1)$\, such that \[\frac{1}{\mu\,\left(\,f(x) \;,\;f(y) \;,\; t\,\right)} - 1\;\leq\;k\left(\,\frac{1}{\mu\,(\,x \;,\;y \;,\; t\,)} - 1\,\right)\] \[and \;\; \frac{1}{\nu\,\left(\,f(x) \;,\;f(y) \;,\; t\,\right)} - 1\;\geq\frac{1}{k}\;\left(\,\frac{1}{\nu\,(\,x \;,\;y \;,\; t\,)} - 1\,\right)\;\;\;\;\;\;\] for each $x \,,\, y\,\in\,X$\, and \,$t \,>\, 0.$\, $($ $k$ is called the contractive constant of\, $f$.$)$
\end{Definition}

%\smallskip

\begin{Proposition}
\cite{Abdul}
Let \,$(\,X \,,\, A\,)$\, be a intuitionistic fuzzy metric space. If \,$f:X\rightarrow\;X$ is fuzzy contractive then f is t-uniformly continuous.
\end{Proposition}

%\smallskip

\begin{Definition}
\cite{Abdul}
Let \,$(\,X \,,\, A\,)$\, be an intuitionistic fuzzy metric space.We will say that the sequence \,$\{\,x_{n}\,\}$\, in \,$X$\, is intuitionistic fuzzy contractive if there \,$k\in\,(0 \,,\, 1)$\, such that \[{\hspace{2.5cm}}\frac{1}{\mu\,(\,x_{n+1} \;,\; x_{n+2} \;,\; t\,)} - 1\;\leq\;k\left(\,\frac{1}{\mu\,(\,x_{n} \;,\; x_{n+1} \;,\; t\,)} - 1\,\right)\] \[and {\hspace{2.5cm}}\frac{1}{\nu\,(\,x_{n+1} \;,\; x_{n+2} \;,\; t\,)} - 1\;\geq\frac{1}{k}\;\left(\,\frac{1}{\nu\,(\,x_{n} \;,\; x_{n+1} \;,\; t\,)} - 1\,\right)\;\;\;\;\] for all \,$t \,>\, 0 ,$\, \,$n\,\in\,N .$
\end{Definition}

%\smallskip

\begin{Theorem}
\cite{Abdul}
Let \,$(\,X \,,\, A\,)$\, be a complete intuitionistic fuzzy metric space with the property \,$\spadesuit$\, in which intuitionistic fuzzy contractive sequences are Cauchy. Let \,$T : X \rightarrow X$\, be a intuitionistic fuzzy contractive mapping such that \,$k$\, is the contractive constant. Then \,$T$\, has a unique fixed point.
\end{Theorem}

%=======================================================================================

\section{ Fixed-point theorems}

%=====================================================================================

\begin{Definition}
Let \,$(\,X \,,\, A\,)$\, be an \textit{\textbf{IFMS}} , $ x\,\in\,X,$\,
$r\,\in\,(0 \,,\, 1)$ , \,$t \,>\, 0$,
 \\${\hspace{2.0cm}}B(\,x \,,\, r \,,\, t\,)\,=\,\left\{\,y\,\in \,X\,/\,\mu(\,x \,,\, y \,,\, t\,) \,>\, 1\,-\,r \,,\, \nu(\,x \,,\, y \,,\, t\,) \,<\, r\,\right\}$.
 \\ Then \,$B(\,x \,,\, r \,,\, t\,)$\, is called an \textbf{open ball} centered at x of radius r w.r.t. t.
\end{Definition}

\begin{Definition}
Let \,$(\,X \,,\, A\,)$\, be an \textit{\textbf{IFMS}} and \,$P\,\subseteq \,X$\,. P is said to be a \textbf{closed set} in \,$(\,X \,,\, A\,)$\, if and only if any sequence \,$\{\,x_{n}\}$\, in P converges to \,$x\,\in P$\, i.e, iff. $\mathop {\lim }\limits_{n\; \to \;\infty } \,\mu\,(\,x_{n} \,,\, x \,,\, t\,)\,=\,1$\, and \,$\mathop {\lim }\limits_{n\; \to \;\infty } \,\nu\,(\,x_{n} \,,\,x \,,\, t\,)\,=\,0\;$  $\Rightarrow$ $x\in\,P$.
\end{Definition}

\begin{Definition}
Let \,$(\,X \,,\, A\,)$\, be an \textit{\textbf{IFMS,}} $ x\,\in\,X,$\,
$r\,\in\,(0 \,,\, 1)$ , \,$t \,>\, 0$,
 \\${\hspace{2.0cm}}S(\,x \,,\, r \,,\, t\,)\,=\,\left\{\,y\,\in \,X\,/\,\mu(\,x \,,\, y \,,\, t\,) \,>\, 1\,-\,r \,,\, \nu(\,x \,,\, y \,,\, t\,) \,<\, r\,\right\}$.
 \\Hence \,$S(\,x \,,\, r \,,\, t\,)$\, is said to be a \textbf{closed ball} centered at x of radius r w.r.t. t iff. any sequence $\{\,x_{n}\}$ in \,$S(\,x \,,\, r \,,\, t\,)$\, converges to y then $y\in\,\,S(\,x \,,\, r \,,\, t\,)$.
\end{Definition}

\begin{Theorem}
$($Contraction on a closed ball$)$ :- Suppose \,$(\,X \,,\, A\,)$\, is a complete IFMS with the property $\spadesuit$ in which IF contractive sequences are Cauchy. Let $T:X\,\rightarrow\,X$ be IF contractive mapping on \,$S(\,x_{\,0} \,,\, r \,,\, t\,)$ \,with contractive constant k. Moreover, assume that
\[{\hspace{3.7cm}}\frac{1}{\mu\,\left(\,x_{0} \,,\,T(x_{0}) \,,\, t\,\right)} \,-\, 1
\,<\, (1 \,-\, k)\,\left(\,\frac{1}{1 \,-\, r}\,-\,1\,\right)\]
\[and {\hspace{2.0cm}}\frac{1}{\nu\,\left(\,x_{0} \,,\,T(x_{0}) \,,\, t\,\right)}
\,-\,1 \,>\, \frac{1}{1 \,-\, k}\,\left(\,\frac{1}{r}-1\,\right)\]
Then T has unique  fixed point in \,$S(\,x_{\,0} \,,\, r \,,\, t\,)$ .
\end{Theorem}

{\bf Proof.}  Let $x_{1}\,=\,T(x_{0}) \,,\, x_{2}\,=\,T(x_{1})\,=\,T^{\,2}(x_{0}) \,,\; \cdots \;,\, x_{n}\,=\,T(x_{n-1})$ \\i.e, \,$\,x_{n}\,=\,T^{\,n}(x_{0})$\, for all $n\,\in\,N.\;$ Now \[{\hspace{2.5cm}}\frac{1}{\mu\,(\,x_{0} \,,\, x_{1} \,,\, t\,)} \,-\, 1 \,<\, (1 \,-\, k)\,\left(\,\frac{1}{1 \,-\,r} \,-\,1\,\right)  \,=\, (1 \,-\, k)\,\left(\,\frac{r}{1 \,-\, r}\,\right)\]
\[i.e,\hspace{0.5cm}\,\frac{1}{\mu\,(\,x_{0} \,,\, x_{1} \,,\, t\,)} \,<\, (1 \,-\, k)\,\left(\,\frac{r}{1 \,-\, r}\,\right) \;+\; 1 {\hspace{2.8cm}}\]
\[ {\hspace{2.0cm}} =\;\frac{r \,-\, r\,k \,+\, 1 \,-\, r}{1 \,-\, r}\,=\,\frac{1 \,-\, r\,k}{1 \,-\, r}\]
\[\Rightarrow\;\;{\mu\,(\,x_{0} \,,\, x_{1} \,,\, t\,)} \,>\, \frac{1 \,-\, r}{1 \,-\, r\,k} \,>\, 1 \,-\, r {\hspace{2.2cm}}\]
\[\Rightarrow\;\;{\mu\,(\,x_{0} \,,\, x_{1} \,,\, t\,)} \,>\, 1 \,-\, r \hspace{1.5cm} \cdots \hspace{1.5cm}(\,i\,)\]
Again,
\[\frac{1}{\nu\,(\,x_{0} \,,\, x_{1} \,,\, t\,)} \,-\, 1 \,>\, \frac{1}{1 \,-\, k}\,\left(\,\frac{1 \,-\, r}{r}\,\right) {\hspace{2.1cm}}\]
\[i.e. \,,\hspace{0.5cm}\frac{1}{\nu\,(\,x_{0} \,,\, x_{1} \,,\, t\,)} \,>\, \frac{1 \,-\, r}{r\,(1 \,-\, k)} \;+\; 1 {\hspace{3.5cm}}\]
\[ {\hspace{3.2cm}}=\,\frac{1 \,-\, r \,+\, r \,-\, r\,k}{r\,(1 \,-\, k)}\,=\,\frac{1 \,-\, r\,k}{r\,(1 \,-\, k)}\]
\[ \Rightarrow\;\;{\nu\,(\,x_{0} \,,\, x_{1} \,,\, t\,)}\,<\,\frac{r\,(1 \,-\, k)}{1 \,-\, r\,k} \,<\, r {\hspace{4.2cm}}\]
\[\Rightarrow\;\;{\nu\,(\,x_{0} \,,\, x_{1} \,,\, t\,)}\,<\,r  \hspace{1.5cm}
\cdots  \hspace{2.5cm} (\,ii\,) {\hspace{1.3cm}}\]
$(\,i\,)$ and $(\,ii\,)$ $\;\Rightarrow\; x_{1}\,\in\,S(\,x_{\,0} \,,\, r \,,\, t\,).$
\\\\Assume that \,$x_{1} \,,\, x_{2} \,,\; \cdots \;,\, x_{n \,-\, 1}\,\in\,S(\,x_{\,0} \,,\, r \,,\, t\,)$.
We show that $x_{n}\,\in\,S(\,x_{\,0} \,,\, r \,,\, t\,).$
\[\frac{1}{\mu\,(\,x_{1} \,,\, x_{2} \,,\, t\,)} \,-\, 1\;=\;\frac{1}{\mu\,\left(\,T(x_{0}) \,,\, T(x_{1}) \,,\, t\,\right)} \,-\, 1 {\hspace{3.2cm}}\]
\[\leq\;k\,\left(\,\frac{1}{\mu\,(\,x_{0} , x_{1} , t\,)} \,-\, 1\,\right)\]
\[ {\hspace{4.2cm}} <\,k\,\left(\frac{1}{1 \,-\, r} \,-\,1\,\right)\,\,[\,since \,\mu\,(\,x_{0} \,,\, x_{1} \,,\, t\,) > 1 \,-\, r \,]\]
\[\Rightarrow\;\frac{1}{\mu\,(\,x_{1} \,,\, x_{2} \,,\, t\,)}\;<\;k\,\left(\,\frac{1 \,-\, 1 \,+\, r}{1 \,-\, r}\,\right) \;+\; 1\;=\;\frac{r\,k}{1 \,-\, r}\;+\;1\;=\;\frac{r\,k \,+\, 1 \,-\, r}{1 \,-\, r}\]
\[i.e.\,,\hspace{0.5cm}\frac{1}{\mu\,(\,x_{1} \,,\, x_{2} \,,\, t\,)}\;<\;\frac{r\,(\,k \,-\, 1\,) \,+\, 1}{1 \,-\, r} {\hspace{6.5cm}}\]
\[\Rightarrow\;\mu\,(\,x_{1} \,,\, x_{2} \,,\, t\,)\;>\;\frac{1 \,-\, r}{1 \,+\, r\,(\,k \,-\, 1\,)}\;=\;\frac{1 \,-\, r}{1 \,-\, r\,(\,1 \,-\, k\,)}\;>\;1 \,-\, r {\hspace{2.1cm}}\]
$\Rightarrow\;\mu\,(\,x_{1} \,,\, x_{2} \,,\, t\,)\;>\;1 \,-\, r$
\\Again,
\[\frac{1}{\nu\,(\,x_{1} \,,\, x_{2} \,,\, t\,)} \,-\, 1\;=\;\frac{1}{\nu\,\left(\,T(x_{0}) \,,\, T(x_{1}) \,,\, t\,\right)} \,-\, 1 {\hspace{3.4cm}}\]
\[ \geq\;\frac{1}{k}\left(\,\frac{1}{\nu\,(\,x_{0} \,,\, x_{1} \,,\, t\,)} \,-\, 1\,\right)\]
\[i.e.\,, \hspace{0.3cm}\frac{1}{\nu\,(\,x_{1} \,,\, x_{2} \,,\, t\,)} \,-\, 1\;\geq\;\frac{1}{k}\,\left(\frac{1}{\nu\,(\,x_{0} \,,\, x_{1} \,,\, t\,)} \,-\, 1\,\right) {\hspace{5.9cm}}\]
\[ {\hspace{0.9cm}}>\;\frac{1}{k}\,\left(\,\frac{1}{r} \,-\, 1\,\right)\;=\;\frac{1}{k}\;\left(\,\frac{1 \,-\, r}{r}\,\right)\]
\[\Rightarrow\;\;\frac{1}{\nu\,(\,x_{1} \,,\, x_{2} \,,\, t\,)} \,>\, \frac{1 \,-\, r}{r\,k} \;+\; 1\;=\;\frac{1 \,-\, r \,+\, r\,k}{r\,k}\;=\;\frac{1 \,-\, r\,(\,1 \,-\, k\,)}{r\,k} {\hspace{1.5cm}}\]
\[i.e.\,,\hspace{0.5cm}\nu\,(\,x_{1} \,,\, x_{2} \,,\, t\,)\;<\;\frac{k\,r}{1 \,-\, r\,(1 \,-\, k)}\;<\;r {\hspace{3.5cm}}\]
\[\Rightarrow\;\; \nu\,(\,x_{1} \,,\, x_{2} \,,\, t\,) \;<\; r {\hspace{6.5cm}}\]
Similarly it can be shown that ,\\
$\mu\,(\,x_{2} \,,\, x_{3} \,,\, t\,) \;>\; 1 \,-\, r \;,\; \nu\,(\,x_{2} \,,\, x_{3} \,,\, t\,) \;<\; r$ \,,\,
$\; \cdots\; ,\, \mu\,(\,x_{n \,-\, 1} \,,\, x_{n} \,,\, t\,) \;>\; 1 \,-\, r$ and
$\;\nu\,(\,x_{n \,-\, 1} \,,\, x_{n} \,,\, t\,) \;<\; r.$ \\Thus, we see that ,
\[\mu\,(\,x_{0} \,,\, x_{n} \,,\, t\,)\;\geq\;\mu\,\left(\,x_{0} \,,\, x_{1} \,,\, \frac{t}{n}\,\right)\;\ast\;\mu\,\left(\,x_{1} \,,\, x_{2} \,,\,
\frac{t}{n}\,\right)\;\ast\;\cdots\;\ast\;\mu\,\left(x_{n \,-\, 1} \,,\, x_{n} \,,\, \frac{t}{n}\,\right)\]
\[\;\;>\; (1 \,-\, r)\;\ast\;(1 \,-\, r)\;\ast\;\cdots\;\ast\;(1 \,-\, r)\;=\; 1 \,-\, r\]
\[i.e. \,, \hspace{0.5cm}\mu\,(\,x_{0} \,,\, x_{n} \,,\, t\,) \;>\; 1 \,-\, r {\hspace{9.5cm}}\]
\[\nu\,(\,x_{0} \,,\, x_{n} \,,\, t\,)\;\leq\;\nu\,\left(\,x_{0} \,,\, x_{1} \,,\, \frac{t}{n}\,\right)\;\diamond\;\nu\,\left(\,x_{1} \,,\, x_{2} \,,\,
\frac{t}{n}\,\right)\;\diamond\;\cdots\;\diamond\;\nu\,\left(x_{n \,-\, 1} \,,\, x_{n} \,,\, \frac{t}{n}\,\right)\]
\[\;<\;r\;\diamond\;r\;\diamond\;\cdots\;\diamond\;r\;=\;r {\hspace{4.2cm}}\]
Thus , $\;\mu\,(\,x_{0} \,,\, x_{n} \,,\, t\,) \;>\; 1 \,-\, r \;and\; \nu\,(\,x_{0} \,,\, x_{n} \,,\, t\,) \;<\; r$
\\$\Rightarrow\;x_{n}\;\in\;S(\,x_{\,0} \,,\, r \,,\, t\,)$
\\Hence, by the theorem 22\cite{Abdul}, T has unique fixed point in \,$S(\,x_{\,0} \,,\, r \,,\, t\,)$.

\begin{Note}
It follows from the proof of Theorem 2.13\cite{Abdul} that for any $\;x\in\;X$ the sequence of iterates \,$\{\,T^{\,n}(x)\}$\, converges to the fixed point of T.
\end{Note}

\begin{Lemma}
Let \,$(\,X \,,\, A\,)$\, be \textit{\textbf{IFMS}} and $T:\,X\rightarrow\,X$ be t-uniformly  continuous on X. If $\;x_{n}\;\rightarrow\;x\; $ as $\; n\;\rightarrow\infty\;$\, in \,$(\,X \,,\, A\,)$\, then $\;T(x_{n})\;\rightarrow\;T(x)\;$ as $\; n\;\rightarrow\infty\;$ in \,$(\,X \,,\, A\,)$ .
\end{Lemma}
{\bf Proof.} \, Proof directly follows from the definitions of t-uniformly continuity and convergence of a sequence in a \textit{\textbf{IFMS}}.

\begin{Lemma}
Let \,$(\,X \,,\, A\,)$\, be \textit{\textbf{IFMS}}. If $\;x_{n}\;\rightarrow\;x$ and $\;y_{n}\;\rightarrow\;y$ in \,$(\,X \,,\, A\,)$\, then $\mu\,(\,x_{n} \,,\, y_{n} \,,\, t\,)\;\rightarrow\;\mu\,(\,x \,,\, y \,,\, t\,) \;and\; \nu\,(\,x_{n} \,,\, y_{n} \,,\, t\,)\;\rightarrow\;\nu\,(\,x \,,\, y \,,\, t\,) \;as \;n\;\rightarrow\infty\; for\; all\; t \,>\, 0\; in \;R\;.$
\end{Lemma}

{\bf Proof.} \,\,We have,
\[ {\hspace{1.0cm}}\mathop {\lim }\limits_{n\;\, \to \;\,\infty } \;\mu\,(\,x_{n} \;,\,x \;,\; t\,)\;=\;1\; , \;\mathop {\lim }\limits_{n\;\, \to \;\,\infty } \;\nu\,(\,x_{n} \;,\,x \;,\; t\,)\;=\;0\;\]
\[ and \;\;\mathop {\lim }\limits_{n\;\, \to \;\,\infty } \;\mu\,(\,y_{n} \;,\,y \;,\; t\,)\;=\;1\; , \;\mathop {\lim }\limits_{n\;\, \to \;\,\infty } \;\nu\,(\,y_{n} \;,\,y \;,\; t\,)\;=\;0\;\]
\[ \mu\,(\,x_{n} \,,\, y_{n} \,,\, t\,)\;\geq\;\mu\,\left(\,x_{n} \,,\, x \,,\, \frac{t}{2}\,\right)\;\ast\;\mu\,\left(\,x \,,\,y_{n} \,,\, \frac{t}{2}\,\right)\]
\[{\hspace{4.8cm}}\;\geq\;\mu\,\left(\,x_{n} \,,\, x \,,\, \frac{t}{2}\,\right)\;\ast\;\mu\,\left(\,x \,,\, y \,,\, \frac{t}{4}\,\right)\;\ast\;\mu\,\left(\,y \,,\, y_{n} \,,\, \frac{t}{4}\,\right)\]
\\${\hspace{0.5cm}}\Rightarrow\;\mathop {\lim }\limits_{n\;\, \to \;\,\infty } \;\mu\,(\,x_{n} \;,\,y_{n} \;,\; t\,)\;\geq\;\mu\,(\,x\;,\;y\;,\;t\,)\;$
\[ \mu\,(\,x \,,\, y \,,\, t\,)\;\geq\;\mu\,\left(\,x \,,\, x_{n} \,,\, \frac{t}{2}\,\right)\;\ast\;\mu\,\left(\,x_{n} \,,\, y \,,\, \frac{t}{2}\,\right) {\hspace{2.5cm}}\]
\[ {\hspace{2.7cm}}\;\geq\;\mu\,\left(\,x \,,\, x_{n} \,,\, \frac{t}{2}\,\right)\;\ast\;\mu\,\left(\,x_{n} \,,\, y_{n} \,,\, \frac{t}{4}\,\right)\;\ast\;\mu\,\left(\,y_{n} \,,\, y \,,\, \frac{t}{4}\,\right)\]
${\hspace{0.5cm}}\Rightarrow\;\mu\,(\,x \;,\,y \;,\; t\,)\;\geq\;\;\mathop {\lim }\limits_{n\;\, \to \;\,\infty } \;\mu\,(\,x_{n}\;,\;y_{n}\;,\;t\,)\;$ $\;\forall\;t\,\,>\,0.$
\\Then,\\ ${\hspace{0.9cm}}\;\mathop {\lim }\limits_{n\;\, \to \;\,\infty } \;\mu\,(\,x_{n} \;,\;y_{n} \;,\; t\,)\;=\;\mu\,(\,x\;,\;y\;,\;t\,)\;$ for all $\;t\;>\,0$,
\\\\Similarly, $\;\mathop {\lim }\limits_{n\;\, \to \;\,\infty } \;\nu\,(\,x_{n} \;,\;y_{n} \;,\; t\,)\;=\;\nu\,(\,x\;,\;y\;,\;t\,)\;$ for all $\;t\;>\,0$.

\begin{Theorem}
\,$(\,X \,,\, A\,)$\, be a complete \textit{\textbf{IFMS}} with the property $\, \spadesuit\,$ in which IF contractive sequences are Cauchy sequences and $T:X\rightarrow X$ be a t-uniformly continuous on X. If for same positive integer m, \,$T^{\,m}\,$ is a IF contractive mapping with k its contractive constant then T has a unique fixed point.
\end{Theorem}

{\bf Proof.}Let $B\;=\;T^{\,m}$,  n be a arbitrary but fixed positive integer and $x_{0}\;\in\;X$.\\we now show that $\;B^{\,n}\,T(x_{0}) \;\rightarrow\; B^{\,n}\,(x_{0})$ in \,$(\,X \,,\, A\,)$\,.
\\Now,
\[\frac{1}{\mu\,\left(\,B^{\,n}T(x_{0}) \;,\; B^{\,n}(x_{0}) \;,\; t\,\right)} \;-\; 1\;=\;\frac{1}{\mu\,\left(\,B(B^{\,n \,-\, 1}\,T(x_{0})) \;,\; B(B^{\,n \,-\, 1}(x_{0})) \;,\; t\,\right)} \;-\; 1\]
\[{\hspace{5.6cm}}\;\leq\;k\,\left(\,\frac{1}{\mu\,(\,B^{n \,-\, 1}\,T(x_{0}) \;,\; B^{n \,-\, 1}(x_{0}) \;,\; t\,)} \;-\; 1\,\right) \;\;\leq\;\;\cdots\;\;\]
\[{\hspace{3.6cm}}\;\leq\;k^{\,n}\,\left(\,\frac{1}{\mu\,(\,T(x_{0}) \;,\; x_{0} \;,\; t\,)} \;-\; 1\,\right)\]
\[\Longrightarrow \;\; 0\;\leq\;\mathop {\lim }\limits_{n\;\, \to \;\,\infty } \left(\,\frac{1}{\mu\,\left(\,B^{\,n}T(x_{0}) \;,\; B^{\,n}(x_{0}) \;,\; t\,\right)} \;-\; 1\,\right)\;\leq\;0 {\hspace{2.9cm}}\]
\[\Longrightarrow\;\;\mathop {\lim }\limits_{n\;\, \to \;\,\infty }\mu\,\left(\,B^{\,n}T(x_{0}) \;,\; B^{\,n}(x_{0}) \;,\; t\,\right)\;=\;1,\;\; for \;all \;\;t\;>\;0. {\hspace{2.5cm}}\]
\\Similarly,
${\hspace{0.5cm}}\mathop {\lim }\limits_{n\;\, \to \;\,\infty }\nu\,\left(\,B^{\,n}\,T(x_{0}) \;,\; B^{\,n}(x_{0}) \;,\; t\,\right)\;=\;0,\;$ for all $\; t\;>\;0$. \\
Thus, $\;B^{\,n}\,T(x_{0}) \;\rightarrow\; B^{\,n}\,(x_{0})$ in \,$(\,X \,,\, A\,)$\,.
\\Again, by the theorem 22\cite{Abdul}, we see that B has a unique fixed point x(say), and from the note [3.5], it follows that $\;B^{\,n}(x_{0})\,\rightarrow \,x$ as $\;n\;\rightarrow\;\infty\;$ in \,$(\,X \,,\, A\,)$\,.
\\\\Since T is t-uniformly continuous on X, it follows from the above lemma[3.6] that $\,B^{\,n}\,T(x_{0})\;=\;T\;B^{\,n}(x_{0})\;\rightarrow\;T(x)$ as $n\;\rightarrow\;\infty$ in \,$(\,X \,,\, A\,)$\,.
\\\\$\Rightarrow\mathop {\lim }\limits_{n\;\, \to \;\,\infty }\mu\,\left(\,B^{\,n}\,T(x_{0}) \;,\; B^{\,n}(x_{0}) \;,\; t\,\right)\;=\;1\;$ and $\;\mathop {\lim }\limits_{n\;\, \to \;\,\infty }\nu\,\left(\,B^{\,n}\,T(x_{0}) \;,\; B^{\,n}(x_{0}) \;,\; t\,\right)\;=\;0$.  By lemma [3.7], we have
\\$\mathop {\lim }\limits_{n\;\, \to \;\,\infty }\mu\,\left(\,T(x) \;,\; x \;,\; t\,\right)\;=\;1$ and $\mathop {\lim }\limits_{n\;\, \to \;\,\infty }\nu\,\left(\,T(x) \;,\; x \;,\; t\,\right)\;=\;0,\;$ for all $\;t\;>\;0$,
\\$i.e. \,, \;\;\mu\,\left(\,T(x) \;,\; x \;,\; t \,\right)\;=\;1\;$ and $\;\nu\,\left(\,T(x) \;,\; x \;,\; t\,\right)\;=\;0,\;$ for all $\;t\;>\;0$.
\\$\Rightarrow\;T(x)\;=\;x\;$
$\;\Rightarrow\;x$ is a fixed point of T. \\If $x^{'}$ is a fixed point of T,
$\,i.e.\,, \;T(x^{'})\;=\;x^{'}$, then $\;T^{\,m}(x^{'})\;=\;T^{\,m \,-\, 1}(T(x^{'}))\;=\;T^{\,m \,-\, 1}(x^{'})\;=\; \cdots \;=\;x^{'}$
$\Rightarrow\;B(x^{'})\;=\;x^{'}\;\Rightarrow\;x^{'}\,$ is a fixed point of B.
 \\But x is the unique fixed point of B, therefore $x\;=\;x^{'}$ which implies that x is the unique fixed point of T.

\begin{Definition}
Let \,$(\,X \,,\, A\,)$\, be \textit{\textbf{IFMS}} and $\;T:X\rightarrow X$.
T is said to be TS-IF contractive mapping if there exists \,$k\;\in\;(0 \,,\, 1)$ such that \[k\;\mu\,\left(\,T(x) \;,\; T(y) \;,\; t\,\right)\;\geq\;\mu\,(\,x \;,\; y \;,\; t\,)\;\] \[ and  \hspace{0.8cm}\frac{1}{k}\;\nu\,\left(\,T(x) \;,\; T(y) \;,\; t\,\right)\;\leq\;\nu\,(\,x \;,\; y \;,\; t\,) \;\;\forall \;\;t \,>\, 0.\]
\end{Definition}

\begin{Theorem}
Let \,$(\,X \,,\, A\,)$\, be a complete \textit{\textbf{IFMS}} with the property $\;\spadesuit\;$ and $\;T : X \rightarrow X$ be TS-IF contractive mapping with k its contraction constant. Then T has a unique fixed point.
\end{Theorem}

{\bf Proof.} Let $x\,\in\,X$ and $x_{\,n}\;=\;T^{\,n}(x)$\, for all \,$n\,\in\,N$. Now for each $t\,>\,0,$
\[k\;\mu\,(\,x_{2} \;,\; x_{1}
\;,\; t\,)\;=\;k\;\mu\,\left(\,T(x_{1}) \;,\; T(x) \;,\; t\,\right) {\hspace{3.5cm}}\]
\[\;\geq\;\mu\,(\,x_{1} \;,\; x \;,\; t\,) {\hspace{2.2cm}}\]
\[i.e. \,,\;\;k\;\mu\,(\,x_{2} \;,\; x_{1} \;,\; t\,)\;\geq\;\mu\,(\,x_{1} \;,\; x \;,\; t\,) {\hspace{6.2cm}}\]
and \[\frac{1}{k}\;\nu\,(\,x_{2} \;,\; x_{1}
\;,\; t\,)\;=\;\frac{1}{k}\;\nu\,\left(\,T(x_{1}) \;,\; T(x) \;,\; t\,\right) {\hspace{3.5cm}}\]
\[\;\leq\;\nu\,(\,x_{1} \;,\; x \;,\; t\,) {\hspace{2.2cm}}\]
\[i.e. \,, \;\;\frac{1}{k}\;\nu\,(\,x_{2} \;,\; x_{1}
\;,\; t\,)\;\leq\;\nu\,(\,x_{1} \;,\; x \;,\; t\,) {\hspace{6.2cm}}\]
Again,\[k\;\mu\,(\,x_{3} \;,\; x_{2} \;,\; t\,)\;=\;k\;\mu\,\left(\,T(x_{2}) \;,\; T(x_{1}) \;,\; t\,\right) {\hspace{3.5cm}}\]
\[\;\geq\;\mu\,(\,x_{2} \;,\; x_{1} \;,\; t\,){\hspace{2.2cm}}\]
\[\Rightarrow\;\; k^{2}\;\mu\,(\,x_{3} \;,\; x_{2} \;,\; t\,)\;\geq\;k\;\mu\,(\,x_{2} \;,\; x_{1} \;,\; t\,)\;\geq\;\;\mu\,(\,x_{1} \;,\; x \;,\;t \,) {\hspace{1.5cm}}\]
\[i.e. \,, \hspace{0.2cm}k^{2}\;\mu\,(\;x_{3} \;,\; x_{2} \;,\; t\,)\;\geq\;\mu\,(\,x_{1} \;,\; x \;,\; t\,) {\hspace{5.5cm}}\]
and \[\frac{1}{k}\;\nu\,(\,x_{3} \;,\; x_{2} \;,\; t\,)\;=\;\frac{1}{k}\;\nu\,\left(\,T(x_{2}) \;,\; T(x_{1}) \;,\; t\,\right) {\hspace{3.3cm}}\]
\[\;\leq\;\nu\,(\,x_{2} \;,\; x_{1} \;,\; t\,) {\hspace{2.0cm}}\]
\[\Rightarrow\;\;\frac{1}{k^{\,2}}\;\nu\,(\,x_{3} \;,\; x_{2} \;,\; t\,)\;\leq\;\frac{1}{k}\nu\,(\,x_{2} \;,\; x_{1} \;,\; t\,)\;\leq\;\nu\,(\,x_{1} \;,\; x \;,\; t\,) {\hspace{2.5cm}}\]
\[i.e. \,, \;\;\frac{1}{k^{\,2}}\;\nu\,(\,x_{3} \;,\; x_{2} \;,\;
t\,)\;\leq\;\nu\,(\,x_{1} \;,\; x \;,\; t\,) {\hspace{6.0cm}}\]
By Mathematical induction, we have,
\[k^{\,n}\;\mu\,(\,x_{\,n \,+\, 1} \;,\; x_{\,n} \;,\; t\,)\;\geq\;\mu\,(\,x_{1} \;,\; x \;,\; t\,)\;\;\;
and \;\;\; \] \[\frac{1}{k^{\,n}}\;\nu\,(\,x_{\,n \,+\, 1} \;,\; x_{\,n} \;,\; t\,)\;\leq\;\nu\,(\,x_{1} \;,\; x \;,\; t\,),\;\; for \;all \; t \,>\, 0.\]
\\We now verify that \,$\{\,x_{\,n}\}$\, is a cauchy sequence in \,$(\,X \,,\, A\,)$. Let
$t_1 \;=\; \frac{t}{p}.$
\\$\mu(\;x_{n},\;x_{n+p},\;t)$
\[\;\geq\;\mu\,(\,x_{n} \;,\;x_{n \,+\, 1} \;,\; t_{1}\,)\;\ast\;\mu(\;x_{n \,+\, 1} \;,\; x_{n \,+\, 2} \;,\; t_{1}\,)\;
\ast\; \cdots \;\ast\;\mu\,(\,x_{n \,+\, p \,-\, 1} \;,\; x_{n \,+\, p} \;,\;t_{1})\]
\[=\;\left(\,\frac{1}{k^{\,n}}\;k^{\,n}\;\mu\,(\,x_{n} \;,\; x_{n \,+\, 1} \;,\; t_{1}\,)\,\right)\;\ast\;
\left(\,\frac{1}{k^{\,n \,+\, 1}}\;k^{\,n \,+\, 1}\;\mu\,(\,x_{n \,+\, 1} \;,\; x_{n \,+\, 2} \;,\; t_{1}\,)\,\right)\;\ast\; \cdots \]
\[{\hspace{6.5cm}}\ast\;\left(\,\frac{1}{k^{\,n \,+\, p \,-\, 1}}\;k^{ \,n \,+\, p \,-\, 1}\;\mu\,(\,x_{\,n \,+\, p \,-\, 1} \;,\; x_{\,n \,+\, p} \;,\; t_{1}\,)\,\right)\]
\[\geq\;\left(\,\frac{1}{k^{\,n}}\;\mu\,(\,x_{1} \;,\; x \;,\; t_1\,)\,\right)\;\ast\;\left(\,\frac{1}{k^{\,n \,+\, 1}}\;
\mu\,(\,x_{1} \;,\; x \;,\; t_1\,)\,\right)\;\ast\; \cdots \;\ast\;\left(\,\frac{1}{k^{\,n \,+\, p \,-\, 1}}\;
\mu\,(\,x_{1} \;,\; x \;,\; t_1\,)\,\right)\]
\[\geq\;\left(\,\frac{1}{k^{\,n}}\;\mu\,(\,x_{1} \;,\; x \;,\; t_1\,)\,\right) \;\ast\; \left(\,\frac{1}{k^{\,n}}\;\mu\,(\,x_{1} \;,\; x \;,\; t_1\,)\,\right) \;=\; \left(\,\frac{1}{k^{\,n}}\;\mu\,(\,x_{1} \;,\; x \;,\; t_1\,)\,\right) \;\;\]
\[\Rightarrow\;1\;<\;\mathop {\lim }\limits_{n\;\, \to \;\,\infty }\;\left(\,\frac{1}{k^{\,n}}\;\mu\,(\,x_{1} \;,\; x \;,\; t_1\,)\,\right) \;\leq\; \mathop {\lim }\limits_{n\;\, \to \;\,\infty }\;\mu\,(\,x_{n} \;,\; x_{n \,+\,p} \;,\; t\,) \;\leq\; 1 {\hspace{1.2cm}}\]
\\$\Rightarrow\;\mathop {\lim }\limits_{n\;\, \to \;\,\infty }\mu\,(\,x_{n} \;,\; x_{n \,+\, p} \;,\; t\,)\;=\;1$
\\\\In the similar way, we have,
\\$\Rightarrow\;\mathop {\lim }\limits_{n\;\, \to \;\,\infty }\nu\,(\,x_{n} \;,\; x_{n \,+\, p} \;,\; t\,)\;=\;0$
\\\\Hence \,$\{\,x_{\,n}\}_{n}$\, is a cauchy sequence in \,$(\,X \,,\, A\,)$.
\\ So, $\;\exists \;y\in\,X$\, such that \,$x_{n}\;\rightarrow\;y$\, as \,$n\;\rightarrow\;\infty$\, in \,$(\,X \,,\, A\,)$.
\\Now,
\[k\;\mu\,\left(\,T(x_{n}) \;,\; T(y) \;,\; t\,\right)\;\geq\;\mu\,(\,x_{n} \;,\; y \;,\; t\,) {\hspace{2.5cm}}\]
\[i.e. \,,\; \;\mu\,\left(\,T(x_{n}) \;,\; T(y) \;,\; t\,\right)\;\geq\;\frac{1}{k}\;\mu\,(\,x_{n} \;,\; y \;,\; t\,) {\hspace{2.6cm}}\]
\[ \Rightarrow\;\;\mathop {\lim }\limits_{n\; \to \;\infty} \mu\,\left(\,T(x_{n}) \;,\; T(y) \;,\; t\,\right)\;\geq\;\mathop {\lim }\limits_{n\; \to \;\infty }\;\frac{1}{k}\;\mu\,(\,x_{n} \;,\; y \;,\; t\,)\;=\;\frac{1}{k}\;>\;1\]
\[\Rightarrow \;\; 1\;<\;\mathop {\lim }\limits_{n\; \to \;\infty} \mu\,\left(\,T(x_{n}) \;,\; T(x) \;,\; t\,\right)\;\leq\;1 {\hspace{5.0cm}}\]
\[\Rightarrow \;\;\mathop {\lim }\limits_{n\; \to \;\infty} \mu\,\left(\,T(x_{n}) \;,\; T(x) \;,\; t\,\right)\;=\;1 {\hspace{5.9cm}}\]
Again,
\[\nu\,\left(\,T(x_{n}) \;,\; T(y) \;,\; t\,\right) \;\leq\; k\;\nu\,(\,x_{n} \;,\; y \;,\; t\,) {\hspace{2.5cm}}\]
\[ \Rightarrow\;\;\mathop {\lim }\limits_{n\; \to \;\infty} \nu\,\left(\,T(x_{n}) \;,\; T(y) \;,\; t\,\right)\;\leq\;\mathop {\lim }\limits_{n\; \to \;\infty }\;k\;\nu\,(\,x_{n} \;,\; y \;,\; t\,)\;=\;0 \;\;\]
\[\Rightarrow \;\;\mathop {\lim }\limits_{n\; \to \;\infty} \nu\,\left(\,T(x_{n})
 \;,\; T(x) \;,\; t\,\right)\;=\;0 {\hspace{5.0cm}}\]
Thus we see that
\[\mathop {\lim }\limits_{n\; \to \;\infty} \mu\,\left(\,T(x_{n}) \,,\,T(y) \,,\, t\,\right)\;=\;1 \;and\; \mathop {\lim }\limits_{n\; \to \;\infty}\nu\,\left(\,T(x_{n}) \,,\, T(y) \,,\, t\,\right)\;=\;0, \;for \; all\; t\,>\,0.\]
\[\Rightarrow\;\mathop {\lim }\limits_{n\; \to \;\infty}\;T(x_{n})\;=\;T(y) \;in \;(\,X \,,\, A\,) \;\;\Longrightarrow \;\;\mathop {\lim }\limits_{n\; \to \;\infty}\;x_{n+1}\;=\;T(y) \;in \;(\,X \,,\, A\,).\]
$i.e. \,, \;y\;=\;T(y)$
\\
$\Rightarrow\;\;$ y is a fixed point of T.
\\To prove the uniqueness, assume
$T(z)\,=\,z\;$ for some $\;z\;\in\;X$. Then  for $t\;>\;0,$ we have
\[\mu\,(\,y \,,\, z \,,\, t\,)\;=\;\mu\,\left(\,T(y) \,,\, T(z) \,,\, t\,\right) {\hspace{2.0cm}}\]
\[\geq\;\frac{1}{k}\;\mu\,(\,y \,,\, z \,,\, t\,) {\hspace{0.6cm}}\]
\[{\hspace{0.5cm}}=\;\frac{1}{k}\;\mu\,\left(\,T(y) \,,\, T(z) \,,\, t\,\right)\]
\[{\hspace{0.8cm}}\geq\;\frac{1}{k^{\,2}}\;\mu\,(\,y \,,\, z \,,\, t\,) \; \geq \;\;\cdots \]
\[ {\hspace{3.0cm}} \geq\;\frac{1}{k^{\,n}}\;\mu\,(\,y \,,\, z \,,\, t\,)\;\longrightarrow\;\infty \,\;as \;\,n\rightarrow\;\infty\]
\[ {\hspace{0.5cm}}\Longrightarrow \;\;1\;<\;\mathop {\lim }\limits_{n\; \to \;\infty}\;\frac{1}{k^{\,n}}\;\mu\,(\,y \,,\, z \,,\, t\,)\;\leq\;\mu\,(\,y \,,\, z \,,\, t\,)\;\leq\;1\]
$\Rightarrow\;\; \mu\,(\,y \,,\, z \,,\, t\,)\;=\;1$
\[\nu\,(\,y \,,\, z \,,\, t\,)\;=\;\nu\,\left(\,T(y) \,,\, T(z) \,,\, t\,\right) {\hspace{2.0cm}}\]
\[\leq\;k\;\mu\,(\,y \,,\, z \,,\, t\,) {\hspace{0.6cm}}\]
\[{\hspace{0.5cm}}=\;k\;\nu\,\left(\,T(y) \,,\, T(z) \,,\, t\,\right)\]
\[{\hspace{0.8cm}}\leq\;k^{\,2}\;\nu\,(\,y \,,\, z \,,\, t\,) \; \leq \;\;\cdots \]
\[ {\hspace{3.0cm}} \leq\;k^{\,n}\;\mu\,(\,y \,,\, z \,,\, t\,)\;\longrightarrow\;0
 \,\;as \;\,n\rightarrow\;\infty\]
\[ {\hspace{0.5cm}}\Longrightarrow \;\;0\;\leq\;\mu\,(\,y \,,\, z \,,\, t\,)
 \;\leq\;\mathop {\lim }\limits_{n\; \to \;\infty}\;k^{\,n}\;\mu\,(\,y \,,\, z \,,\, t\,)\;<\;0\]
$\Rightarrow\;\;\nu\,(\,y \,,\, z \,,\, t\,)\;=\;0$
\\$\Rightarrow\;\;y\;=\;z$
\\This completes the proof.

\end{document}